\documentclass[letterpaper, twoside, 11pt]{amsart}

\usepackage{amsmath,amsthm,amsfonts,amssymb,amscd}
\usepackage{bbm}
\usepackage{tikz}
\usepackage{tikz-cd}
\usepackage{standalone}
\usepackage[hidelinks]{hyperref}
\usetikzlibrary{arrows,chains,matrix,positioning,scopes}
\usepackage[letterpaper, left=3cm,right=3cm, top=3cm, bottom=3cm]{geometry}
\usepackage{color}

\newcommand{\QQ}{{\mathbb Q}}

\newcommand{\fu}{{\mathfrak{u}}}

\newcommand\dual{\raise0.9ex\hbox{$\scriptscriptstyle\vee$}}

\theoremstyle{definition}

\theoremstyle{plain}
\newtheorem{thm}{Theorem}

\newtheorem{lemma}[thm]{Lemma}  
\newtheorem{cor}[thm]{Corollary}
\numberwithin{thm}{section}

\theoremstyle{remark}

\makeatletter
\def\@seccntformat#1{%
  \protect\textup{\protect\@secnumfont
    \ifnum\pdfstrcmp{subsection}{#1}=0 \bfseries\fi% subsection # in \bfseries
    \csname the#1\endcsname
    \protect\@secnumpunct
  }%
}  
\makeatother

\tikzset{>=stealth}

\begin{document}

\title{On endomorphisms of extensions in Tannakian categories}
\author{Payman Eskandari}
\address{Department of Mathematics and Statistics, University of Winnipeg, Winnipeg MB, Canada }
\email{p.eskandari@uwinnipeg.ca}

\begin{abstract}
We prove some analogues of Schur's lemma for endomorphisms of extensions in Tannakian categories. More precisely, let $\mathbf{T}$ be a neutral Tannakian category over a field of characteristic zero. Let $E$ be an extension of $A$ by $B$ in $\mathbf{T}$. We consider conditions under which every endomorphism of $E$ that stabilizes $B$ induces a scalar map on $A\oplus B$. We give a result in this direction in the general setting of arbitrary $\mathbf{T}$ and $E$, and then a stronger result when $\mathbf{T}$ is filtered and the associated graded objects to $A$ and $B$ satisfy some conditions. We also discuss the sharpness of the results.
\end{abstract}
\maketitle

\section{Introduction}
Let $K$ be a field of characteristic zero and $\mathbf{T}$ a Tannakian\footnote{Throughout the paper, all Tannakian categories are neutral. We will freely use the language of Tannakian categories. See \cite{DM} for a reference.} category over $K$. Given any object $X$ of $\mathbf{T}$, let $End_\mathbf{T}(X)$ be the endomorphism algebra of $X$. Given a subobject $Y$ of $X$, denote the subalgebra of $End_\mathbf{T}(X)$ consisting of the endomorphisms that restrict to an endomorphism of $Y$ (i.e. that map $Y$ to $Y$) by $End_\mathbf{T}(X;Y)$. 

Let $A$ and $B$ be nonzero objects of $\mathbf{T}$. Fix an extension of $A$ by $B$: 
\begin{equation}\label{eq0}
\begin{tikzcd}
0 \arrow[r] & B \arrow[r] & E \arrow[r] & A \arrow[r] & 0.
\end{tikzcd}
\end{equation}
In this note we prove some analogues of Schur's lemma for $End_\mathbf{T}(E;B)$. 

The extension \eqref{eq0} induces a homomorphism of algebras
\begin{equation}\label{eq1}
\Omega: End_\mathbf{T}(E;B) \rightarrow End_\mathbf{T}(B)\times End_\mathbf{T}(A) \hspace{0.5in} \phi\mapsto(\phi_B,\phi_A),
\end{equation}
where given $\phi\in End_\mathbf{T}(E;B)$, its image $(\phi_B,\phi_A)$ is characterized by the commutativity of
\begin{equation}\label{eq2}
\begin{tikzcd}
  0 \arrow[r] & B  \arrow[d, "\phi_{B}"] \arrow[r, ] & E \arrow[d, "\phi"] \arrow[r, ] & A \arrow[d, "\phi_A"] \arrow[r] & 0 \\
    0 \arrow[r] & B  \arrow[r, ] & E \arrow[r, ] &  A  \arrow[r] & 0 .
\end{tikzcd}
\end{equation}
The image of $\Omega$ always contains the diagonal copy of $K$ in $End_\mathbf{T}(B)\times End_\mathbf{T}(A)$ (as the image of scalar endomorphisms of $E$). Roughly speaking, it is natural to expect that the further away \eqref{eq0} is from splitting, the smaller the image of $\Omega$ should be. We shall prove two results in this spirit. The first is the following:
\begin{thm}\label{thm1}
Let $\fu(E)$ be the Lie algebra of the kernel of the homomorphism from the Tannakian group of $E$ to the Tannakian group of $A\oplus B$, naturally considered as a subobject of the internal Hom object $Hom(A,B)$ (see \S \ref{sec: review of results on Tannakian groups of extensions} below for details). Assume that $\fu(E)=Hom(A,B)$. Then the image of $\Omega$ is equal to the diagonal copy of $K$.
\end{thm}
The works \cite{Ha06} and \cite{Ha11} of Hardouin (in the case where $A$ and $B$ are semisimple) and \cite{EM1} of the author and Kumar Murty (for arbitrary possibly non-semisimple $A$ and $B$) give a characterization of the subobject $\fu(E)$ of $Hom(A,B)$. A summary of this characterization is recalled in \S \ref{sec: review of results on Tannakian groups of extensions} below. It follows from this characterization that the condition that $\fu(E)=Hom(A,B)$, which we refer to as the maximality of $\fu(E)$, implies that the extension class
\[\mathcal{E} \in Ext^1_\mathbf{T}(\mathbbm{1}, Hom(A,B))\]
(where $Ext^1_\mathbf{T}$ is the $Ext^1$ group in $\mathbf{T}$ and $\mathbbm{1}$ is the unit object) corresponding to \eqref{eq0} under the canonical isomorphism
\[
Ext^1_\mathbf{T}(A,B) \cong Ext^1_\mathbf{T}(\mathbbm{1}, Hom(A,B))
\]
is {\it totally nonsplit}, i.e. for any proper subobject $C$ of $Hom(A,B)$ the pushforward of $\mathcal{E}$ along the quotient $Hom(A,B)\rightarrow Hom(A,B)/C$ is nonsplit.

In the case where $A$ and $B$ are semisimple, the maximality of $\fu(E)$ is equivalent to the total nonsplitting of $\mathcal{E}$. But in general, the two conditions are not equivalent, as the examples in \S \ref{further remarks} illustrate. The second result of the paper asserts that in some important settings one can relax the hypothesis of Theorem \ref{thm1} from assuming maximality of $\fu(E)$ to assuming total nonsplitting of $\mathcal{E}$:

%Let $\mathbf{T}^{ss}$ be the full subcategory of all the semisimple objects of $\mathbf{T}$. Fix\footnote{Choosing such a functor amounts to choosing a Levi factor of the Tannakian group of $\mathcal{G}$ with respect to a fiber functor. It can always be done since $K$ has characteristic zero.} a
%an exact faithful ($K$-linear) tensor functor $\mathbf{T}\rightarrow \mathbf{T}^{ss}$ restricting to the identity on $\mathbf{T}^{ss}$. Denote the image of an object $X$ under the functor by $X^{ss}$. The second result of the paper is the following:
%
%\begin{thm}\label{thm3}
%Suppose that the extension $\mathcal{E}$ (defined as above) is totally nonsplit. Suppose moreover that there are no nonzero morphisms between $\mathbf{1}$ and $\inHom(A,B)^{ss}$. Then the image of $\Omega$ is equal to the diagonal copy of $K$.
%\end{thm}

%Note that choosing the functor amounts to choosing a Levi factor of the Tannakian group of $\mathbf{T}$.

\begin{thm}\label{thm2}
Suppose that $\mathbf{T}$ is a filtered Tannakian category with the filtration denoted by $W_\bullet$. Suppose moreover that condition (i) or (ii) below holds:
\begin{itemize}
\item[(i)] The associated graded $Gr^W E$ is semisimple and there are no nonzero morphisms $Gr^W A\rightarrow Gr^W B$.
%$\mathbbm{1}\rightarrow Hom(Gr^W A,Gr^W B)$.
\item[(ii)] The sets of weights of $A$ and $B$ are disjoint. 
\end{itemize} 
Then if $\mathcal{E}$ (defined as above) is totally nonsplit, the image of $\Omega$ will be equal to the diagonal copy of $K$.
%Suppose that $\mathbf{T}$ is a filtered Tannakian category. Suppose moreover that the sets of weights of $A$ and $B$ are disjoint. Assume that $\mathcal{E}$ (defined as above) is totally nonsplit. Then the image of $\Omega$ is equal to the diagonal copy of $K$.
\end{thm}

Here, by a filtered Tannakian category we mean a Tannakian category equipped with a filtration $W_\bullet$ similar to the weight filtration on mixed motives (indexed by integers, functorial, increasing, finite on every object, exact, and respecting the tensor structure). The filtration is referred to as the weight filtration. By the weights of an object $X$ we mean the integers $n$ such that $W_nX/W_{n-1}X$ is not zero. The associated graded of $X$ is defined to be $Gr^W X:=\bigoplus_n W_nX/W_{n-1}X$. The prototype examples of filtered Tannakian categories are various Tannakian categories of mixed motives and the category of mixed Hodge structures. In the former categories always (assuming the category is a reasonable category of mixed motives) and in the latter category if $E$ is graded polarizable, then $Gr^WE$ is semisimple. Of course, it is only useful to include condition (ii) as a separate condition in the statement if $Gr^WE$ is not known to be semisimple.

The kernel of $\Omega$ is canonically isomorphic to $Hom_\mathbf{T}(A,B)$, where $Hom_\mathbf{T}$ is the Hom group in $\mathbf{T}$. Since the functor that sends an object $X$ to $Gr^WX$ is faithful, under condition (i) or (ii) of Theorem \ref{thm2} $Hom_\mathbf{T}(A,B)$ will be zero. Thus Theorem \ref{thm2} in fact implies the following statement\footnote{In the original statement of Theorem \ref{thm2} we wrote the assertion as merely $Im(\Omega)=K$ rather than $End_\mathbf{T}(E;B)\cong K$ because the proof of the result also applies in some situations where the former assertion might actually be weaker than the latter. See \S \ref{other situations where thm 2 holds}, Theorem \ref{thm3}.}:

\begin{cor}
Suppose that $\mathbf{T}$ is a filtered Tannakian category and that condition (i) or (ii) of Theorem \ref{thm2} hold. If $\mathcal{E}$ is totally nonsplit, then 
\[
End_\mathbf{T}(E;B)\cong K
\]
and every element of $End_\mathbf{T}(E;B)$ is a scalar endomorphisms of $E$.
\end{cor}

Below, we first recall the characterization of $\fu(E)$ mentioned above, and then prove Theorems \ref{thm1} and \ref{thm2}. The final section of the note includes some further remarks. In particular, we give an example that shows that in the general setting of Theorem \ref{thm1} one cannot relax the maximality condition to total nonsplitting. Also, we discuss an example involving 1-motives that shows that in the setting where $\mathbf{T}$ is filtered and the sets of weights of $A$ and $B$ are disjoint, the total nonsplitting of $\mathcal{E}$ does not imply maximality of $\fu(E)$, so that in this setting the second theorem is indeed stronger than the first one. We also discuss a generalization of Theorem \ref{thm2} (see \S \ref{other situations where thm 2 holds}).

{\bf Acknowledgements.} I would like to thank Kumar Murty for many helpful discussions. 

\section{Recollections on Tannakian groups of extensions}\label{sec: review of results on Tannakian groups of extensions}
To simplify the notation, we fix a choice of fiber functor and identify $\mathbf{T}$ with the category of finite dimensional (algebraic) representations of an affine group scheme $\mathcal{G}$ over $K$ (with $\mathcal{G}$ = the Tannakian group of $\mathbf{T}$ with respect to the fiber functor). We will use the same symbol for an object of $\mathbf{T}$ and its underlying vector space\footnote{Or rather, its image under the fiber functor if one doesn't like to identify $\mathbf{T}$ with the category of representations of $\mathcal{G}$.}. For any object $X$ of $\mathbf{T}$ and any $g\in \mathcal{G}$, we denote the image of $g$ in $GL(X)$ by $g_X$. The image of $\mathcal{G}$ in $GL(X)$ is denoted by $\mathcal{G}(X)$; this is the Tannakian group of the Tannakian subcategory $\langle X\rangle$ generated by $X$.\footnote{Recall that this means $\langle X\rangle $ is the smallest full Tannakian subcategory of $\mathbf{T}$ which contains $X$ and is closed under taking subquotients.}

We should point out that even though we think of $\mathbf{T}$ as the category of representations of $\mathcal{G}$, all the objects in $\mathbf{T}$ that appear in the following text (in particular, the object $\fu(E)$ introduced below) will be intrinsic to the Tannakian category $\mathbf{T}$. For this reason, we often prefer to use the terms {\it object} and {\it subobject} ( = object and subobject in $\mathbf{T}$) instead of the terms {\it $\mathcal{G}$-representation} and {\it $\mathcal{G}$-subrepresentation}.

As they were in the Introduction, the Ext and Hom groups in $\mathbf{T}$ are denoted by $Ext_\mathbf{T}$ and $Hom_\mathbf{T}$. We use the notations $Hom$ and $End$ (without any decorations) to refer to the Hom and End groups in the category of vector spaces. As we have adopted the convention of denoting an object  of $\mathbf{T}$ and its underlying vector space by the same symbol, for any objects $X$ and $Y$ of $\mathbf{T}$ the notation $Hom(X,Y)$ will refer to both the internal Hom (which is an object of $\mathbf{T}$) and the Hom space in the category of vector spaces between the underlying vector spaces. This should not lead to confusion as the relevant interpretation will be clear from the context.

Given a vector space $X$ and a subspace $Y$ of $X$, denote the subalgebra of $End(X)$ consisting of linear maps on $X$ which map $Y$ to $Y$ by $End(X;Y)$. Similarly, the subgroup of $GL(X)$ consisting of the elements which map $Y$ to itself is denoted by $GL(X;Y)$. Given an object $X$ of any category, the identity map on $X$ is denoted by $Id_X$. We will sometimes simply write $Id$ if $X$ is clear from the context. 

Fix objects $A$, $B$ and $E$ of $\mathbf{T}$ and the exact sequence \eqref{eq0}. Let $\mathcal{U}(E)$ be the kernel of the natural map
\begin{equation}\label{eq3}
\mathcal{G}(E) \ \twoheadrightarrow \ \mathcal{G}(B\oplus A).
\end{equation}
Choosing a section of $E\twoheadrightarrow B$ in the category of vector spaces to identify 
\[
E = B\oplus A
\]
as vector spaces, we have an embedding
\[
\mathcal{U}(E) \ \rightarrow \ W_{-1}GL(B\oplus A;B):=\{\begin{pmatrix}
Id_B & f \\
0 &Id_A\end{pmatrix}: f\in Hom(A,B)\}.
\]
The group $W_{-1}GL(B\oplus A;B)$ is unipotent and abelian and hence so is $\mathcal{U}(E)$. Since \linebreak $W_{-1}GL(B\oplus A;B)$ is abelian, the embedding above is actually canonical, i.e. does not depend on the choice of the section of $E \twoheadrightarrow A$ used to identify $E=B\oplus A$.%Indeed, a different choice of section would result in conjugating the first embedding by an element of $W_{-1}GL(B\oplus A;B)$, which is an abelian group.

Let $\fu(E)$ be the Lie algebra of $\mathcal{U}(E)$. Then $\fu(E)$ can be identified as a subspace of $Hom(A,B)$ with the exponential map $\fu(E) \rightarrow \mathcal{U}(E)$ simply sending
\[
f\in \fu(E)\subset Hom(A,B)
\hspace{.3in} \text{to} \hspace{.3in}
\begin{pmatrix}
Id & f \\
0 &Id\end{pmatrix}.
\]
Through the adjoint representation of $\mathcal{G}(E)$ on $\fu(E)$, the Lie algebra $\fu(E)$ is naturally equipped with a $\mathcal{G}$-action. The inclusion $\fu(E)\subset Hom(A,B)$ is compatible with the $\mathcal{G}$-actions, making $\fu(E)$ a subobject of the internal Hom $Hom(A,B)$ (see \cite[\S 3.1]{EM1}, for instance). This subobject has a nice description, which we recall now. 

As in  \S 1, let
\[\mathcal{E} \in Ext^1_\mathbf{T}(\mathbbm{1}, Hom(A,B))\]
be the element corresponding to the class of \eqref{eq0} under the canonical isomorphism
\[
Ext^1_\mathbf{T}(A,B) \cong Ext^1_\mathbf{T}(\mathbbm{1}, Hom(A,B)).
\]
For any subobject $C\subset Hom(A,B)$, the pushforward of $\mathcal{E}$ along the quotient map 
\[Hom(A,B)\rightarrow Hom(A,B)/C\]
is denoted by $\mathcal{E}/C$. The following characterization of $\fu(E)$ was proved in \cite{EM1}:

\begin{thm}[Theorem 3.3.1 of \cite{EM1}, see also Lemma 3.4.3 of \cite{EM2}]\label{thm with Kumar} 
Given any subobject $C$ of $Hom(A,B)$, we have $\fu(E)\subset C$ if and only if 
the pushforward 
\[
\mathcal{E}/C \in Ext^1_\mathbf{T}(\mathbbm{1}, Hom(A,B)/C)
\]
is in the image of the natural injection
\[
Ext^1_{\langle A\oplus B\rangle }(\mathbbm{1}, Hom(A,B)/C) \, \rightarrow \, Ext^1_\mathbf{T}(\mathbbm{1}, Hom(A,B)/C),
\]
where $Ext^1_{\langle A\oplus B\rangle }$ is the $Ext^1$ group in the Tannakian subcategory $\langle A\oplus B\rangle $ of $\mathbf{T}$ generated by $A\oplus B$. (Thus  $\fu(E)$ is the smallest subobject of $Hom(A,B)$ with this property.)
\end{thm}

In the case where $A$ and $B$ are semisimple, this was earlier proved by Bertrand in \cite{Be01} in the setting of D-modules, and by Hardouin in \cite{Ha06} and \cite{Ha11} in the setting of arbitrary Tannakian categories. In this case, the statement simplifies to the following:

\begin{thm}[Theorem 2 of \cite{Ha11}]\label{Hardouin's thm}
Suppose $A$ and $B$ are semisimple. Let $\mathcal{E}$ be as above. Then given any subobject $C$ of $Hom(A,B)$, we have $\fu(E)\subset C$ if and only if 
the pushforward $\mathcal{E}/C$ splits.
\end{thm}
Note that in the general case (where $A$ or $B$ is not semisimple), by Theorem \ref{thm with Kumar} if $C$ is any subobject of $Hom(A,B)$ such that $\mathcal{E}/C$ splits, then $C$ contains $\fu(E)$. The pushforward $\mathcal{E}/\fu(E)$ however may not split. See the examples in \S \ref{further remarks} below.

We also recall an explicit description of $\mathcal{E}$ (see \cite[\S 3.2]{EM1} for details). Let 
%\begin{align*}
%Hom(A,E)^\dagger & :=\\
%& \hspace{-.3in} \{f\in Hom(A,E): \text{the composition $A\xrightarrow{f}E\twoheadrightarrow A$ is a scalar multiple of $Id_A$}\}.
%\end{align*}
\[
Hom(A,E)^\dagger:= \{f\in Hom(A,E): \text{the composition $A\xrightarrow{f}E\twoheadrightarrow A$ is a scalar multiple of $Id_A$}\}.
\]
It is easy to see that this is a subobject of $Hom(A,E)$. The element $\mathcal{E}$ is the class of the extension
\begin{equation}\label{eq5}
\begin{tikzcd}
0 \arrow[r] & Hom(A,B) \arrow[r] & Hom(A,E)^\dagger \arrow[r] & \mathbbm{1} \arrow[r] & 0.
\end{tikzcd}
\end{equation}
Here, the injective map is simply the obvious embedding, sending $f\in Hom(A,B)$ to 
\[
A\xrightarrow{f} B \hookrightarrow E.
\]
The surjective map in \eqref{eq5} is the map that sends $f\in Hom(A,E)^\dagger$ to $a\in K$, where 
\[
A\xrightarrow{f}E\twoheadrightarrow A
\]
is $a\cdot Id_A$.

\section{Proofs of Theorems \ref{thm1} and \ref{thm2} for $A=\mathbbm{1}$}

%\subsection{The case $A=\mathbbm{1}$}
The goal of this section is to prove Theorems \ref{thm1} and \ref{thm2} in the case where $A=\mathbbm{1}$; the general case will be deduced from this in the next section. In this case, identifying $Hom(\mathbbm{1},B) =B$ the extension $\mathcal{E}$ is simply given by \eqref{eq0}. Theorem \ref{thm with Kumar} asserts that $\fu(E)$
is the smallest subobject of $B$ such that $\mathcal{E}/\fu(E)$ is the pushforward of an extension of $\mathbbm{1}$ by $B/\fu(E)$ in the subcategory $\langle B\rangle$. If $B$ is semisimple, $\fu(E)$ is the smallest subobject of $B$ such that $B/\fu(E)$ splits.
%In this case, we have
%\[\fu(E)\subset \inHom(\mathbbm{1},B) =B\]
%is the smallest subrepresenation of $B$ such that $\mathcal{E}/\fu(E)$ is the pushforward of an extension of $\mathbbm{1}$ by $B/\fu(E)$ in the subcategory $\langle B\rangle$. If $B$ is semisimple, $\fu(E)$ is the smallest subrepresentation such that $B/\fu(E)$ splits. The extension $\mathcal{E}$ is simply given by \eqref{eq0}.

We first establish a lemma:
\begin{lemma}\label{lem 1}
Assume $A=\mathbbm{1}$. Let $\lambda: E\rightarrow B_0$ be a morphism from $E$ to an object $B_0$, such that $B_0$ belongs to the subcategory $\langle B\rangle$. Then $\fu(E)\subset B\cap \ker(\lambda)$.
%Assume $A=\mathbbm{1}$. Suppose that there exists a morphism $\lambda: E\rightarrow B_0$ for some object $B_0$ in the subcategory $\langle B\rangle$. Then $\fu(E)\subset B\cap \ker(\lambda)$.
\end{lemma}

\begin{proof}
Set $B':=B\cap \ker(\lambda)$. Consider the commutative diagram
\begin{equation}\label{eq4}
\begin{tikzcd}
   & 0 \arrow{d} & 0 \arrow{d} & 0 \arrow{d} &\\
   0 \arrow[r] & B' \arrow{d} \arrow[r, ] & \ker(\lambda)  \arrow[d, ] \arrow[r, ] &  \ker(\lambda)/B'  \arrow{d} \arrow[r] & 0 \\
   0 \arrow[r] & B \arrow{d} \arrow[r, ] & E \arrow{d} \arrow[r, ] &  \mathbbm{1} \arrow{d}  \arrow[r] & 0 \, \\
      0 \arrow[r] & B/B' \arrow{d} \arrow[r, ] & E/\ker(\lambda) \arrow{d} \arrow[r, ] &  E/(B+\ker(\lambda))  \arrow{d}  \arrow[r] & 0 \, ,\\
   & 0& 0 & 0 & 
\end{tikzcd}
\end{equation}
where the maps are inclusions and quotient maps. The rows and columns are exact. 

\underline{Case I}: Suppose $\ker(\lambda)\not\subset B$, so that $B'$ is a proper subobject of $\ker(\lambda)$. Being a nonzero subobject of the unit object, $\ker(\lambda)/B'$ must be isomorphism to $\mathbbm{1}$. Thus $\mathcal{E}$ (= the second row) is the pushforward of an extension of $\mathbbm{1}$ by $B'$ (the first row). It follows that $\mathcal{E}/B'$ splits and $\fu(E)\subset B'$ by Theorem \ref{thm with Kumar}.

\underline{Case II}: Suppose $\ker(\lambda)\subset B$, so that $B'=\ker(\lambda)$. Then the third row of the diagram is the pushforward $\mathcal{E}/B'$. By assumption, $E/\ker(\lambda)$ is in the subcategory generated by $B$. Again $\fu(E)\subset B'$ by Theorem \ref{thm with Kumar}.
\end{proof}

We can now establish Theorems \ref{thm1} and \ref{thm2} in the case that $A=\mathbbm{1}$. Let $\phi\in End_\mathbf{T}(E;B)$. Then $\phi_\mathbbm{1}$ (= the induced map on $\mathbbm{1}$ by $\phi$) is equal to $a\cdot Id_\mathbbm{1}$ for some $a\in K$. The endomorphism $\lambda:=\phi-a \cdot Id_E$ of $E$ then factors through $B$, i.e is the composition with the inclusion $B\hookrightarrow E$ of a morphism $E\rightarrow B$, which we also denote by $\lambda$. 

To obtain Theorem \ref{thm1}, apply the previous lemma to $\lambda$. We get $\fu(E)\subset B\cap \ker(\lambda)$. The assumption that $\fu(E)=B$ thus gives $B\subset \ker(\lambda)$, i.e. $\phi=a\cdot Id$ on $B$, as desired.

We now turn our attention to Theorem \ref{thm2}. We thus further assume that $\mathbf{T}$ is a filtered Tannakian category and that either (i) $Gr^WE$ is semisimple and 
\[
Hom_\mathbf{T}(\mathbbm{1}, Gr^WB) = 0,
\]
or (ii) 0 is not a weight of $B$ (note that $A=\mathbbm{1}$ is pure of weight 0). Both conditions guarantee that the kernel of $\lambda: E\rightarrow B$ cannot be contained in $B$. Indeed, this is simply by weight considerations if (ii) holds. On the other hand, if (i) holds, after applying the associated graded functor the sequence \eqref{eq0} splits. Choosing a section for the sequence (which will be unique because $Hom_\mathbf{T}(\mathbbm{1}, Gr^WB)$ vanishes), if $\ker(\lambda)\subset B$ we have a diagram
\[
\begin{tikzcd}
 &&\mathbbm{1} \arrow[d, hookrightarrow] & \\
  0 \arrow[r] & Gr^W\ker(\lambda) \arrow[d, hookrightarrow]   \arrow[r, ] & Gr^W E  \arrow[r, "\lambda"] & Gr^WB \\
    & Gr^WB  & &
\end{tikzcd}
\]
(with obvious maps and the row being exact). We thus get a nonzero morphism $\mathbbm{1}\rightarrow Gr^WB$. 

Thus $B':=B\cap \ker(\lambda)$  is a proper subobject of $\ker(\lambda)$. Considering the diagram \eqref{eq4} with $B'$ and $\lambda:E\rightarrow B$ as in here, the extension $\mathcal{E}$ is the pushforward of the top row along the inclusion $B'\hookrightarrow B$, so that $\mathcal{E}/B'$ splits. If $\mathcal{E}$ is totally nonsplit, we get $B'=B$ and $B\subset \ker(\lambda)$. Thus we have also established Theorem \ref{thm2} when $A=\mathbbm{1}$.

\section{Proofs of Theorems \ref{thm1} and \ref{thm2} for arbitrary $A$}

%\subsection{The case of an arbitrary $A$}
We now assume that $A$ is arbitrary. The extension $\mathcal{E}$ is now given by \eqref{eq5}. Assume the hypotheses of Theorem \ref{thm1} or \ref{thm2} for the extension given in \eqref{eq0}. Then the hypotheses also hold for the extension given in \eqref{eq5}, i.e. if \eqref{eq5} is taken as our original \eqref{eq0}: To see this for Theorem \ref{thm1}, note that in view of Theorem \ref{thm with Kumar} we have $\fu(E)\subset \fu(Hom(A,E)^\dagger)$, as the subcategory $\langle Hom(A,B)\rangle $ is contained in $\langle A\oplus B\rangle$; to see it for Theorem \ref{thm2} note that 
\[
Hom_\mathbf{T}(\mathbbm{1}, Gr^W Hom(A,B))=Hom_\mathbf{T}(\mathbbm{1}, Hom(Gr^WA,Gr^WB))=Hom_\mathbf{T}(Gr^WA,Gr^WB).
\]
Thus by the special case of the results already proved we know that the image of the map
\begin{equation}\label{eq6}
End_\mathbf{T}(Hom(A,E)^\dagger;Hom(A,B)) \longrightarrow End_\mathbf{T}(Hom(A,B))\times End_\mathbf{T}(\mathbbm{1})
\end{equation}
induced by \eqref{eq5} is the diagonal copy of $K$. Hence the general case of the results will be established if we prove the following lemma:

\begin{lemma}\label{lem reduction to A=1 case}
Suppose the image of \eqref{eq6} is the diagonal copy of $K$. Then so is the image of \eqref{eq1}.
\end{lemma}
\begin{proof}
Let $\phi\in End_\mathbf{T}(E;B)$. We will show that $(\phi_B,\phi_A)$ is in the diagonal copy of $K$. Adding a suitable scalar multiple of $Id_E$ to $\phi$ if necessary, we may assume that $\phi$ is an automorphism (recall that $K$ is of characteristic zero). Let $\phi^\dagger\in End(Hom(A,E))$ be the map that sends any $f\in Hom(A,E)$ to the composition
\[
A\xrightarrow{\phi_A^{-1}} A \xrightarrow{~f~} E \xrightarrow{~\phi~} E.
\]
Since $\phi_A$ and $\phi$ are morphisms in $\mathbf{T}$, so is $\phi^\dagger$. Since $B$ is stable under $\phi$, the map $\phi^\dagger$ stabilizes $Hom(A,B)$. Moreover, if $f$ is in $Hom(A,E)^\dagger$ with $f\pmod{A}=Id_A$, then we have a commutative diagram
\[
\begin{tikzcd}
   A \arrow[r, "\phi_A^{-1}"] & A  \arrow[dr, "Id", swap] \arrow[r, "f"] & E  \arrow[d, tail, twoheadrightarrow] \arrow[r, "\phi"] &  E  \arrow[d, tail, twoheadrightarrow]  \\
 & & A \arrow[r, "\phi_A"] & A \, , 
\end{tikzcd}
\]
so that $\phi^\dagger(f)$ is also in $Hom(A,E)^\dagger$ with $\phi^\dagger(f) \pmod{A} $ being the identity map on $A$. We conclude that:
\begin{itemize}
\item[(i)] $\phi^\dagger$ restricts to an element of $End_\mathbf{T}(Hom(A,E)^\dagger;Hom(A,B))$, and
\item[(ii)] denoting this restriction also by $\phi^\dagger$, we have $\phi^\dagger_\mathbbm{1}=Id$ (where $\phi^\dagger_\mathbbm{1}$ is the map induced on $\mathbbm{1}$ by $\phi^\dagger\in End_\mathbf{T}(Hom(A,E)^\dagger)$). 
\end{itemize}
Since the image of \eqref{eq6} is the diagonal copy of $K$, it follows that the restriction of $\phi^\dagger$ to $Hom(A,B)$ is also the identity map. That is, for every linear map $f: A\rightarrow B$, we have
\[
\phi_Bf\phi_A^{-1} = f.
\]
Since $A$ and $B$ are nonzero, by elementary linear algebra $\phi_A$ and $\phi_B$ are both scalar maps and they are given by multiplication with the same element of $K$.
\end{proof}

\section{Further remarks}\label{further remarks}

\subsection{} If \eqref{eq0} is an arbitrary extension in a general Tannakian category $\mathbf{T}$ (with no extra assumptions on \eqref{eq0} or $\mathbf{T}$), total nonsplitting of $\mathcal{E}$ (= the extension of $\mathbbm{1}$ by $Hom(A,B)$ corresponding to \eqref{eq0} under the canonical isomorphism) does not guarantee that the image of $\Omega$ is $K$. Thus the hypothesis of maximality of $\fu(E)$ in Theorem \ref{thm1} cannot be relaxed to total nonsplitting. 

For example, given any field $K$ of characteristic zero, take $\mathbf{T}$ be the category of finite dimensional algebraic representations of the subgroup $\mathcal{G}$ of $GL_3$ (over $K$) consisting of all the matrices of the form
\[
\begin{pmatrix}
1&a&b\\
&1&a\\
&&1
\end{pmatrix},
\]
where the missing entries are zero. Let $B$ be $K^2$ with the action of $\mathcal{G}$ given by left multiplication by the top left $2\times 2$ submatrix, and $E$ be $K^3$ with the canonical action of $\mathcal{G}$ through left multiplication. We have an embedding  $B\hookrightarrow E$ given by $(x_1,x_2)\mapsto (x_1,x_2,0)$, fitting into a short exact sequence
\[
\begin{tikzcd}
0 \arrow[r] & B \arrow[r] & E \arrow[r] & \mathbbm{1} \arrow[r] & 0,
\end{tikzcd}
\]
with the map $E\twoheadrightarrow \mathbbm{1}$ being projection onto the third coordinate. It is easy to see that the extension above is totally nonsplit. However, $E$ has an endomorphism
\[
\phi: (x_1,x_2,x_3)\mapsto (x_2,x_3,0) 
\]
which stabilizes $B$ but its restriction to $B$ is not a scalar multiple of the identity. 

It is worth mentioning that here, by Theorem \ref{thm1} $\fu(E)$ is not maximal, so that this also gives an example that shows that in general, total nonsplitting of $\mathcal{E}$ does not imply that $\fu(E)$ must be maximal (in particular, in general $\mathcal{E}/\fu(E)$ does not have to split). See the next subsection for a more interesting example that also illustrates this.

\subsection{} Assume that $\mathbf{T}$ is filtered and that $A$ and $B$ have disjoint sets of weights. Then total nonsplitting of $\mathcal{E}$ still does not imply maximality of $\fu(E)$, so that Theorem \ref{thm2} is indeed stronger than Theorem \ref{thm1} in this setting. The example provided in \S 6.10 of \cite{EM2} using the work \cite{JR87} of Jacquinot and Ribet on deficient points on semiabelian varieties illustrates this. If we take $\mathbf{T}$ to be the category of mixed Hodge structures, $E$ to be the 1-motive denoted by $M$ in \S 6.10 of \cite{EM2}, and we take $B=W_{-1}M$ and $A=M/W_{-1}M=\mathbbm{1}$, then the sequence \eqref{eq0} (given by the natural inclusion and quotient maps) is totally nonsplit, the weights of $A$ and $B$ are disjoint, and $\fu(E)$ (which is the same as $\fu_{-1}(M)$ in \S 6.10 of \cite{EM2}) is not maximal. In fact, we have $\fu(E)=0$. See {\it loc. cit.}

Continuing to work in the category of mixed Hodge structures, here we include a somewhat simpler example which avoids using deficient points. Let $J$ be a simple complex abelian variety of positive dimension. Let $N$ be a nonsplit extension of $\mathbbm{1}$ by $H_1(J)$. Then $N\dual(1)$ is an extension of $H^1(J)(1)$ by $\QQ(1)$, which after a choice of polarization can be thought of as an extension of $H_1(J)$ by $\QQ(1)$. Since\footnote{See \cite[Lemma 9.3.8]{SGA7} or \cite[Lemma 6.4.1]{EM2}.} the $Ext^2$ groups vanish in the category of mixed Hodge structures (see \cite{Be83}, for example), there is an object $E$ fitting into the diagram
\[
\begin{tikzcd}
   &  & 0 \arrow{d} & 0 \arrow{d} &\\
   0 \arrow[r] & \QQ(1) \ar[equal]{d} \arrow[r, ] & N\dual(1) \arrow[d, ] \arrow[r, ] &  H_1(J) \arrow{d} \arrow[r] & 0 \\
   0 \arrow[r] & \QQ(1) \arrow[r, ] & E \arrow{d} \arrow[r, ] &  N \arrow{d}  \arrow[r] & 0 \, ,\\
      &  & \mathbbm{1} \arrow{d} \ar[equal]{r} &  \mathbbm{1}  \arrow{d}  &  \\
   & &  0 & 0 & 
\end{tikzcd}
\]
in which the rows and columns are exact and the top row and the right column are nonsplit.

Take the first vertical extension of the diagram to play the role of our \eqref{eq0}; it will also be our $\mathcal{E}$. Then $\mathcal{E}$ is totally nonsplit, as $\QQ(1)$ is the unique maximal proper subobject of $N\dual(1)$ and the pushforward $\mathcal{E}/\QQ(1)$ (= the right column) is nonsplit. On the other hand, $\mathcal{E}/\QQ(1)$ is an extension in the subcategory generated by $N\dual(1)$, hence by Theorem \ref{thm with Kumar} we have $\fu(E)\subset \QQ(1)$. In particular, $\fu(E)$ is not maximal.

% Should I say something about E being a 1-motive?

\subsection{}\label{other situations where thm 2 holds}
In the proof of the case $A=\mathbbm{1}$ of Theorem \ref{thm2} the only place where the filtration on $\mathbf{T}$ and condition (i) or (ii) played a part is when we concluded that the kernel of $\lambda:E\rightarrow B$ (with $\lambda$ as in the proof) is not contained in $B$. Combining with Lemma \ref{lem reduction to A=1 case} we obtain the following generalization of Theorem \ref{thm2}:

\begin{thm}\label{thm3}
Let \eqref{eq0} be an extension in any Tannakian category $\mathbf{T}$ over a field of characteristic 0. Suppose that the kernel of any morphism
\begin{equation}\label{eq7}
Hom(A,E)^\dagger \rightarrow Hom(A,B)
\end{equation}
is not contained in $Hom(A,B)$. Then if $\mathcal{E}$ (i.e. the extension of $\mathbbm{1}$ by $Hom(A,B)$ corresponding to \eqref{eq0}, as before) is totally nonsplit, the image of $\Omega$ will be the diagonal copy of $K$.
\end{thm}

In particular, this can be applied in the following situation: Suppose $\mathcal{R}$ is a reductive subgroup of the group $\mathcal{G}(E)$ ( = the Tannakian group of $\langle E\rangle$). Every object of $\langle E\rangle$ can also be considered as an $\mathcal{R}$-representation. In the (semisimple) category of $\mathcal{R}$-representations, we can choose a splitting of $\mathcal{E}$ to decompose
\[
Hom(A,E)^\dagger \simeq Hom(A,B)\oplus \mathbbm{1}.
\]
Suppose that there are no nonzero $\mathcal{R}$-equivariant maps $A\rightarrow B$, or equivalently
\[
\mathbbm{1} \rightarrow Hom(A,B).
\]
Then the kernel of any morphism \eqref{eq7} in $\mathbf{T}$ cannot be contained in $Hom(A,B)$, and hence the image of $\Omega$ will be the diagonal copy of $K$. In fact, since $Hom_\mathbf{T}(A,B)$ is zero, we get $End_\mathbf{T}(E;B)\cong K$.

Note that this scenario directly generalizes the situation of Theorem \ref{thm2}: If $\mathbf{T}$ is filtered, taking $\mathcal{R}$ to be $\mathcal{G}(Gr^WE)$ embedded in $\mathcal{G}(E)$ via the section of $\mathcal{G}(E)\twoheadrightarrow \mathcal{G}(Gr^WE)$ induced by $Gr^W: \langle E\rangle \rightarrow \langle Gr^WE\rangle$ we recover case (i) of Theorem \ref{thm2}. Taking $\mathcal{R}$ to be the multiplicative group $\mathbb{G}_m$ embedded in $\mathcal{G}(E)$ through a (possibly noncentral) cocharacter $\mathbb{G}_m \rightarrow \mathcal{G}(E)$ that induces the weight grading on the associated gradeds we recover case (ii) of the result.

\subsection{} We have
\[
\ker(\Omega) = Hom_\mathbf{T}(A,B),
\]
where $Hom_\mathbf{T}(A,B)$ is considered as a subset of $End_\mathbf{T}(E)$ via
\[
(A\xrightarrow{~f~}B) \ \ \mapsto \ \ (E\twoheadrightarrow A \xrightarrow{~f~}B \hookrightarrow E).
\]
Whenever $Im(\Omega)=K$, the natural embedding of $K$ into $End_\mathbf{T}(E;B)$ as the space of scalar maps provides a section for the short exact sequence
\[\begin{tikzcd}
0 \arrow[r] & Hom_\mathbf{T}(A,B) \arrow[r] & End_\mathbf{T}(E;B) \arrow[r, "\Omega"] & Im(\Omega) \arrow[r] & 0.
\end{tikzcd}\]
This gives an isomorphism
\[
End_\mathbf{T}(E;B) \cong K\oplus Hom_\mathbf{T}(A,B).
\]
The isomorphism is initially of vector spaces only, but transferring the multiplication on $End_\mathbf{T}(E;B)$ to the right hand side it becomes an isomorphism of algebras. The multiplication on the right is given by 
\[
(a,f)(a',f') = (aa', af'+a'f)
\] 
and the embedding of $K$ is through the first factor. In particular, $End_\mathbf{T}(E;B)$ is commutative if $Im(\Omega)=K$.


\begin{thebibliography}{00}

%\bibitem{An19}
%Y. Andr\'{e},
%A note on 1-motives,
%International Mathematics Research Notices 2021, no. 3, 2074-2080
%
%\bibitem{Andre}
%Y. Andr\'{e},
%Mumford-Tate groups of mixed Hodge structures and the theorem of the fixed part,
%Compositio Mathematica,  Volume 82 (1992) no. 1,  p. 1-24

%\bibitem{ACGH}
%E. Arbarello , M. Cornalba , P. A. Griffiths , J. Harris,
%{\it Geometry of Algebraic Curves}, Volume I,  Grundlehren der mathematischen Wissenschaften Volume 267, Springer

\bibitem{Be83}
A. A. Beilinson, 
Notes on absolute Hodge cohomology, 
Applications of Algebraic K-theory to Algebraic Geometry and Number Theory, Part I, Proceedings of a Summer Research Conference held June 12-18, 1983, in Boulder, Colorado, Contemporary Mathematics 55, American Mathematical Society, Providence, Rhode Island, pp. 35–68

%\bibitem{Be02}
%C. Bertolin,
%The Mumford-Tate group of 1-motives,
%Ann. Inst. Fourier, Grenoble, 52, 4 (2002), 1041-1059
%
%\bibitem{Be03}
%C. Bertolin,
%Le radical unipotent du groupe de Galois motovique d'un 1-motif,
%Math. Ann. 327, 585-607 (2003)

\bibitem{Be01}
D. Bertrand,
Unipotent radicals of differential Galois group and integrals of solutions of inhomogeneous equations,
Math. Ann. 321 (2001), no. 3, 645–666


%\bibitem{Carlson}
%J. A. Carlson,
%Extensions of mixed Hodge structures,
%Journ\'{e}es de G\'{e}ometrie Alg\'{e}brique d'Angers, Juillet 1979/Algebraic Geometry, Angers, 1979, pp. 107–127, Sijthoff \& Noordhoff, Alphen aan den Rijn, Md., 1980.
%
%
%\bibitem{De74}
%P. Deligne,
%Theorie de Hodge III,
%Publications Math\'{e}matiques de l'I.H.\'{E}.S., tome 44 (1974), p. 5-77
%
%\bibitem{De79}
%P. Deligne,
%Valeurs de fonctions L et p\'{e}riodes d’int\'{e}grales,
%Proceedings of Symposia in Pure Mathematics (AMS), Vol. 33 (1979), part 2, pp. 313–346

%\bibitem{De82}
%P. Deligne (Notes by J.S. Milne),
%Hodge cycles on abelian varieties,
%In Hodge Cycles, Motives, and Shimura Varieties,
%Lecture Notes in Mathematics 900, Springer-Verlog, Berlin (1982)
%
\bibitem{DM}
P. Deligne and J.S. Milne,
Tannakian Categories,
In Hodge Cycles, Motives, and Shimura Varieties,
Lecture Notes in Mathematics 900, Springer-Verlog, Berlin (1982)
%
%
%\bibitem{Elk90}
%R. Elkik,
%Le th\'{e}or\`{e}me de Manin-Drinfeld,
%Ast\'{e}risque No. 183 (1990), 59–67.

\bibitem{EM1}
P. Eskandari, V. K. Murty,
The fundamental group of an extension in a Tannakian category and the unipotent radical of the Mumford-Tate group of an open curve,
to appear in the Pacific Journal of Mathematics

\bibitem{EM2}
P. Eskandari, V. K. Murty,
On unipotent radicals of motivic Galois groups,
Algebra \& Number Theory 17-1 (2023), pp 165-215 %DOI 10.2140/ant.2023.17.165


%
%
\bibitem{SGA7}
A. Grothendieck, 
Mod\`{e}les de N\'{e}ron et monodromie, 
SGA VII.1, no 9, Springer LN 288, 1968

\bibitem{Ha06}
C. Hardouin,
Hypertranscendance et Groupes de Galois aux diff\'{e}rences, arXiv 0609646v2, 2006

\bibitem{Ha11}
C. Hardouin,
Unipotent radicals of Tannakian Galois groups in positive characteristic, 
Arithmetic and Galois theories of differential equations, 223-239, S\'{e}min. Congr., 23, Soc. Math. France, Paris, 2011
%

\bibitem{JR87}
O. Jacquinot and K. Ribet, 
Deficient points on extensions of abelian varieties by G, 
J. Number Theory,
25: 2 (1986), 133-151

%\bibitem{Jan88}
%U. Jannsen,
%Deligne cohomology, Hodge-D-conjecture, and motives,
%in {\it Beilinson's Conjectures on Special Values of L-Functions}, Perspectives in Mathematics 4, Academic Press 1988, pp 305-372


%
%\bibitem{Jan90}
%U. Jannsen,
%Mixed motives and algebraic K-theory,
%Lecture Notes in Mathematics 1400, Springer-Verlag, Berlin, 1990


%\bibitem{Jo14}
%P. Jossen,
%On the Mumford-Tate conjecture for 1-motives,
%Inventiones Math. (2014) 195: 393-439



%\bibitem{Voisin}
%C. Voisin,
%{\it Hodge Theory and Complex Algebraic Geometry}, Volume 1,
%Cambridge Studies in Advanced Mathematics 76, 2002


\end{thebibliography}
\end{document}